# A New Algorithm for Linear Programming

Dhananjay P. Mehendale
Sir Parashurambhau College, Tilak Road, Pune-411030,
India

## Abstract

In this paper we propose a new polynomial-time algorithm for linear programming. We further extend the ideas used in this new linear programming algorithm for nonlinear programming problems. The standard linear programming problem consists of an objective function to be optimized, either maximized or minimized, subject to certain constraints provided in terms of inequalities. The new algorithm is based on the idea of treating the objective function as a parameter. We form a matrix of coefficients, made-up of the coefficients of the variables defined in the problem itself, and the coefficients of variables defined newly, for converting inequalities into equations, namely, slack variables if it is the maximization problem, or, surplus variables if it is the minimization problem. The system of equations we use consist of the objective equation and equations obtained from inequalities defining constraint imposed by the problem. We obtain reduced-row-echelon-form, R, for this matrix containing only one variable, namely, the objective function itself as an unknown parameter, d, say. This matrix in the reduced-row- echelon-form contains columns (column vectors) corresponding to basic variables and non-basic variables. If all the entries in the columns corresponding to non-basic variables in R are already nonnegative then we will see that we have almost reached to the solution and nothing much is left to be done. If there are columns corresponding to non-basic variables which contain some negative entries then we will require to apply suitable row transformations, at most $m$ in number if there are m rows in R, as we will see below, to make all the entries in the columns corresponding to non-basic variables nonnegative. We then proceed to show that the method developed above for linear programming naturally extends to nonlinear programming problems. For nonlinear programming problems we use the technique of Grobner bases, since Grobner basis is an equivalent of reduced row echelon form for a system of nonlinear equations.

**1. Introduction:** There are two types of linear programs (linear programming problems):

1. Maximize: $C^T x$
   Subject to: $Ax \leq b$
   $x \geq 0$
   Or
2. Minimize: $C^T x$
   Subject to: $Ax \geq b$
   $x \geq 0$

where $x$ is a column vector of size n×1 of unknown variables. We call these variables the problem variables where $C$ is a column vector of size n×1 of profit (for maximization problem) or cost (for minimization problem) coefficients, and $C^T$ is a row vector of size 1×n obtained by matrix transposition of $C$.

where $A$ is a matrix of constraints coefficients of size m×n.

where $b$ is a column vector of constants of size m×1 representing the boundaries of constraints.

By introducing the appropriate slack variables (for maximization problem) and surplus variables (for minimization problem), the above mentioned linear programs gets converted into standard form as:

$$\begin{aligned} \text{Maximize: } & C^T x \\ \text{Subject to: } & Ax + s = b \\ & x \geq 0, s \geq 0 \end{aligned} \qquad (1.1)$$

where s is slack variable vector of size m×1.
This is a maximization problem.
Or

$$\begin{aligned} \text{Minimize: } & C^T x \\ \text{Subject to: } & Ax - s = b \\ & x \geq 0, s \geq 0 \end{aligned} \qquad (1.2)$$

where $s$ is surplus variable vector of size m×1.
This is a minimization problem.

In geometrical language, the constraints defined by the inequalities form a region in the form of a convex polyhedron, a region bounded by the constraint planes, $Ax_i = b_i$, and the coordinate planes. This region is called feasible region and it is straightforward to establish that there exists at least one vertex of this polyhedron at which the optimal solution for the problem is situated when the problem at hand is well defined, i.e. neither inconsistent, nor unbounded, nor infeasible. There may be unique optimal solution and sometimes there may be infinitely many optimal solutions, e.g. when one of the constraint planes is parallel to the objective plane we may have a multitude of optimal solutions. The points on an entire plane or an entire edge can constitute the optimal solution set.

These problems are handled most popularly by using the well known simplex algorithm or some of its variant. Despite its theoretical exponential complexity the simplex method works quite efficiently for most of the practical problems. However, there are few computational difficulties associated with simplex algorithm. In order to view them in nutshell we begin with stating some common notions and definitions that are prevalent in the literature. A variable $x_i$ is called basic variable in a given equation if it appears with unit coefficient in that equation and with zero coefficients in all other equations. A variable which is not basic is called non-basic variable. A sequence of elementary row operations that changes a given system of linear equations into an equivalent system (having the same solution set) and in which a given non-basic variable can be made a basic variable is called a pivot operation. An equivalent system containing basic and non-basic variables obtained by application of suitable elementary row operations is called canonical system. At times, the introduction of slack variables for obtaining standard form automatically produces a canonical system, containing at least one basic variable in each equation. Sometimes a sequence of pivot operations is needed to be performed to get a canonical system. The solution obtained from canonical system by setting the non-basic variables to zero and solving for the basic variables is called basic solution and in addition when all the variables have nonnegative values the solution satisfying all the imposed constraints is called a basic feasible solution. Simplex method cannot start without an initial basic feasible solution. The process of finding such a solution, which is a necessity in many of practical problems, is called Phase I of the simplex algorithm. Simplex method starts its Phase II with an initial basic feasible solution in canonical form at hand. Then simplex tests whether this solution is optimal by checking whether all the values of relative profits (profits that result due to unit change in the values of non-basic variables) of all the non-basic variables are nonpositive. When not optimal, the simplex method obtains an adjacent basic feasible solution by selecting a non-basic variable having largest relative profit to become basic. Simplex then determines and carries out the exiting of a basic variable, by the so called minimum ratio rule, to change it into a non-basic variable leading to formation of a new canonical system. On this new canonical system the whole procedure is repeated till one arrives at an optimal solution.

The main computational difficulties of the simplex method which may cause the reduction in its computational efficiency are as follows:

1] There can be more than one non-basic variable with largest value for relative profit and so a tie can take place while selecting a non-basic variable to become basic. The choice at this situation is done arbitrarily and so the choice made at this stage causing largest possible per unit improvement is not necessarily the one that gives largest total improvement in the value of the objective function and so not necessarily minimizes the number of simplex iterations.

2] While applying minimum ratio rule it is possible for more than one constraint to give the same least ratio causing a tie in the selection of a basic variable to leave for becoming non-basic. This degeneracy can cause a further complication, namely, the simplex method can go on without any improvement in the objective function and the method may trap into an infinite loop and fail to produce the desired optimal solution. This phenomenon is called cycling which enforces modification in the algorithm by introducing some additional time consuming rules that reduce the efficiency of the simplex algorithm.

3] Simplex is not efficient on theoretical grounds basically because it searches adjacent basic feasible solutions only and all other simplex variants which examine nonadjacent solutions as well have not shown any appreciable change in the overall efficiency of these modified simplex algorithms over the original algorithm.

Because of the far great practical importance of the linear programs and other similar problems in the operations research it is a most desired thing to have an algorithm which works in a single step, if not, in as few steps as possible. No method has been found which will yield an optimal solution to a linear program in a single step ([1], Page 19). We aim to propose a polynomial-time algorithm for linear programming which aims at fulfilling this requirement in a best possible and novel way.

2. **A New Polynomial-time Algorithm for Linear Programming:** In this section we propose and discuss our new algorithm for solving linear programming problems. This algorithm is based on treating objective function as a parameter and aims at finding its optimal value.

We start with the following equation:

$$C^T x = d \qquad (2.1)$$

where $d$ is an unknown parameter, and call it the objective equation. The (parametric) plane defined by this equation will be called objective plane.

Please note that we are discussing first the maximization problems. A similar approach for minimization problems will be discussed next.

Given a maximization problem**,** we first construct the combined system of equations containing the objective equation and the equations defined by the constraints imposed by the problem under consideration, combined into a single matrix equation, viz.,

$$\begin{bmatrix} C^T_{(1\times n)} & 0_{(1\times m)} \\ A_{(m\times n)} & I_{(m\times m)} \end{bmatrix} \begin{bmatrix} x \\ s \end{bmatrix} = \begin{bmatrix} d \\ b \end{bmatrix} \qquad (2.2)$$

Let $\mathbf{E} = \begin{bmatrix} C^T_{(1\times n)} & 0_{(1\times m)} \\ A_{(m\times n)} & I_{(m\times m)} \end{bmatrix}$**,** and let $\mathbf{F} = \begin{bmatrix} d \\ b \end{bmatrix}$

Let **[E, F]** denote the augmented matrix obtained by appending the column vector **F** to matrix **E** as a last column. We then find *R*, the reduced-row-echelon-form (rref) ([2], pages 73-75) of the above augmented matrix **[E, F]**. Thus,

$$R = \text{rref}\,([\mathbf{E}, \mathbf{F}]) \qquad (2.3)$$

Note that the augmented matrix **[E, F]** as well as its reduced row echelon form *R* contains **only one parameter**, namely, $d$ and all other entries are **constants**. The maximization problem of linear programming asks to determine the unique maximum value for $d$, for which there exists a feasible solution and this unique maximal value of the objective function has been used to obtain the feasible solution and so this solution is in fact the optimal solution. In the case of an unbounded linear programming problem there is no upper bound (lower bound, in the case of minimization problem) for the value of $d$, while in the case of an infeasible linear program the set of feasible solutions is empty. The steps that will be executed to determine the optimal solution will also tell by implication when such optimal solution does not exist in the case of an unbounded or infeasible problem.

The general form of the matrix *R* representing the reduced row echelon form is

$$R = \begin{bmatrix} a_{11} & a_{12} & \cdots & a_{1n} & b_{11} & b_{12} & \cdots & b_{1m} & c_1 d + e_1 \\ a_{21} & a_{22} & \cdots & a_{2n} & b_{21} & b_{22} & \cdots & b_{2m} & c_2 d + e_2 \\ a_{31} & a_{32} & \cdots & a_{3n} & b_{31} & b_{32} & \cdots & b_{3m} & c_3 d + e_3 \\ \vdots & \vdots & \vdots & \vdots & \vdots & \vdots & \vdots & \vdots & \vdots \\ a_{m1} & a_{m2} & \cdots & a_{mn} & b_{m1} & b_{m2} & \cdots & b_{mm} & c_m d + e_m \\ a_{(m+1)1} & a_{(m+1)2} & \cdots & a_{(m+1)n} & b_{(m+1)1} & b_{(m+1)2} & \cdots & b_{(m+1)n} & c_{(m+1)} d + e_{(m+1)} \end{bmatrix}$$

The first $n$ columns of the above matrix represent the coefficients of the problem variables (i.e. variables defined in the linear program) $x_1, x_2, \cdots, x_n$. The next $m$ columns represent the coefficients of the slack variables $s_1, s_2, \cdots, s_m$ used to convert inequalities into equalities to obtain the standard form of the linear programming problem. The last column represents the transformed right hand side of the equation (2.2) during the process (a suitable sequence of transformations) that is carried out to obtain the reduced-row-echelon-form. Note that the last column of *R* contains the linear form $d$ as a parameter whose optimal value is to be determined such that the nonnegativity constraints remain valid, i.e. $x_i \geq 0, 1 \leq i \leq n$ and $s_j \geq 0, 1 \leq j \leq m$. Among first $(n+m)$ columns of *R* the columns correspond to basic (pivot) variables are the columns that are unit vectors and the remaining ones to non-basic variables.

Now, if there are some negative entries in the columns of *R* corresponding to non-basic variables then we carry out suitable elementary row transformations on the obtained *R* = rref (**[E, F]**) so that every column among the columns associated with non-basic variables become nonnegative. We are doing this because as will be seen below we can then put zero value for these non-basic variables and can determine the values of all the basic variables and the linear programming problem will then be solved completely. It is easy to check that for a linear programming problem if all the coefficients of parameter $d$ in the last column of *R* are positive then the linear programming problem at hand is unbounded since in such case the parameter $d$ can be increased arbitrarily without violating the nonnegativity constraints on variables $x_i, s_j$. Also, for a linear programming problem if all the coefficients of some non-basic variable represented by a column of *R* are nonpositive and are strictly negative in those rows having a negative coefficient to parameter $d$ that appears in the last column of these rows then again the problem belongs to the category of unbounded problems since we can increase the value of $d$ to any high value without violating the nonnegativity constraints for the variables by assigning sufficiently high value to this non-basic slack variable. Note that the rows of *R* actually offer expressions for basic variables in terms of non-basic variables and terms of type $c_k d + e_k, k = 1,2,\cdots(m+1)$ containing the parameter $d$ on the right side. The rows with a positive coefficient for the parameter $d$ represent those equations in which the parameter $d$ can be increased arbitrarily without violating the nonnegativity constraints on variables $x_i, s_j$. So, these equations with a positive coefficient for the parameter $d$ are not implying any upper bound on the maximum possible value of parameter $d$. However, these rows are useful in certain situations as they are useful to find lower bound on the value of parameter $d$. The rows with a negative coefficient for the parameter $d$ represent those equations in which the parameter $d$ cannot be increased arbitrarily without violating the nonnegativity constraints on variables $x_i, s_j$. So, these equations with a negative coefficient for the parameter $d$ are implying an upper bound on the maximum possible value of parameter $d$ and so important ones for maximization problems. Note that actually every row of *R* is offering us a value for parameter $d$ which can be obtained by equating to zero each term of the type $c_k d + e_k, k = 1,2,\cdots(m+1)$. We now define the sub-matrix of $R$, say $R_N$, made up of all columns of *R* and containing those rows $j$ of *R* for which the coefficients $c_j$ of the parameter $d$ are negative. Let $c_{i_1}, c_{i_2}, \cdots, c_{i_k}$ coefficient of $d$ in the rows of $R$ which are negative. We collect these rows with negative coefficient for $d$ to form the mentioned submatrix, $R_N$, of $R$ given below. With this it is clear that coefficients of $d$ in all other rows of *R* are greater than or equal to zero.

$$R_N = \begin{bmatrix} a_{i_1 1} & a_{i_1 2} & \cdots & a_{i_1 n} & b_{i_1 1} & b_{i_1 2} & \cdots & b_{i_1 m} & c_{i_1} d + e_{i_1} \\ a_{i_2 1} & a_{i_2 2} & \cdots & a_{i_2 n} & b_{i_2 1} & b_{i_2 2} & \cdots & b_{i_2 m} & c_{i_2} d + e_{i_2} \\ a_{i_3 1} & a_{i_3 2} & \cdots & a_{i_3 n} & b_{i_3 1} & b_{i_3 2} & \cdots & b_{i_3 m} & c_{i_3} d + e_{i_3} \\ \vdots & \vdots & \vdots & \vdots & \vdots & \vdots & \vdots & \vdots & \vdots \\ \vdots & \vdots & \cdots & \vdots & \vdots & \vdots & \cdots & \vdots & \vdots \\ a_{i_k 1} & a_{i_k 2} & \cdots & a_{i_k n} & b_{i_k 1} & b_{i_k 2} & \cdots & b_{i_k m} & c_{i_k} d + e_{i_k} \end{bmatrix}$$

It should be clear to see that if $R_N$ is empty (i.e. not containing a single row) then the problem at hand is unbounded. Among the first $(n+m)$ columns of $R_N$ first $n$ columns represent the coefficients of problem variables and next $m$ columns represent the coefficients of slack variables. There are certain columns corresponding to basic (pivot) variables which are unit vectors. The other columns correspond to non-basic variables. As mentioned, among the columns for non-basic variables those having all entries nonnegative can only lead to decrement in the value of $d$ when a positive value is assigned to them. This is undesirable as we aim maximization of the value of $d$. So, we can safely set the value of such variables equal to zero. When all columns corresponding to non-basic variables in $R$ are having all entries nonnegative then we can set all non-basic variables to zero, set $d = \min\{d\}$ in every row of $R$ and find the basic feasible solution which will be optimal, with $\min\{d\}$ as optimal value for the objective function at hand.

In $R_N$ we now proceed to consider those non-basic variables for which the columns of $R_N$ contain some (at least one) positive values and some negative (at least one) values. In such case when we assign some positive value to such non-basic variable it leads to decrease in the value of $d$ in those rows in which $c_k > 0$ and increase in the value of $d$ in those rows in which $c_k < 0$. We now need to consider the ways of dealing with this situation. We deal with this situation as follows: In this case, we choose and carry out appropriate and legal elementary row transformations on the matrix $R$ in the reduced row echelon form to achieve nonnegative value for all the entries in the columns corresponding to non-basic variables in the matrix $R$. The elementary row transformations are chosen to produce new matrix which remains equivalent to original matrix in the sense that the solution set of the matrix equation with original matrix and matrix equation with transformed matrix remain same. Due to this equivalence we can now set all the non-basic variables in this transformed matrix to zero and obtain with justification $d_{\min} = \min\{d\}$ as optimal value for the objective function and obtain basic feasible solution as optimal solution by substitution.

Let us now discuss our new algorithm in steps:

**Algorithm 2.1 (Maximization):**

1. Express the given problem in standard form:

   Maximize: $C^T x$

   Subject to: $Ax + s = b$

   $x \geq 0, s \geq 0$

2. Construct the augmented matrix [**E F**], where

   $$\mathbf{E} = \begin{bmatrix} C^T_{(1\times n)} & 0_{(1\times m)} \\ A_{(m\times n)} & I_{(m\times m)} \end{bmatrix}, \text{ and } \mathbf{F} = \begin{bmatrix} d \\ b \end{bmatrix}$$

   and obtain the reduced row echelon form:

   $R$ = rref ([**E, F**])

3. If there is a row (or rows) of zeroes at the bottom of $R$ in the first n columns and containing a nonzero constant in the last column then declare that the problem is inconsistent and stop.

4. **4.** Else if for any value of $d$ one observes that nonnegativity constraint for some variable gets violated for at least one of the variables then declare that the problem at hand is infeasible and stop.
5. **5.** Else, if the coefficients of $d$ in the last column are all positive or if there exists a column of R corresponding to some non-basic variable with all entries negative then declare that the problem at hand is unbounded and stop.
6. **6.** Else check whether all the entries in the columns of $R$ corresponding to non-basic variables are nonnegative.
7. **7.** If yes, then find the sub-matrix of *R*, say $R_N$, made up of those rows of *R* for which the coefficient of $d$ in the last column is negative and solve $c_{i_r} d + e_{i_r} = 0$ for each such a term in the last column of $R_N$ and find the value of $d = d_{i_r}$ for $r = 1,2,\cdots,k$ and find $d_{\min} = \min\{ d_{i_r} \}$ which represents the maximal value of the objective function.
8. **8.** Set the value for all non-basic variables equal to zero.
9. **9.** Substitute the value $d_{\min} = \min\{ d_{i_r} \}$ in place of $d = d_{i_r}$ in every row of R and find the values for all basic (pivot) variables. The set of values assigned to basic as well as non-basic variables together forms the optimal solution for the problem.
10. **10.** If not, then apply suitable elementary row transformations (which do exist when the problem at hand is a well-defined legitimate problem) first on $R_N$ and make all the entries in the columns corresponding to non-basic variables nonnegative. Now, apply suitable elementary row transformations on $R$ if still there are some negative entries in the column corresponding to non-basic variables and change $R$ such that now every column corresponding to non-basic variables will strictly contain nonnegative entries in $R$. Now perform steps 7-9 in the new $R_N$ contained in the new $R$, which now contains nonnegative entries in every column corresponding to non-basic variable, and find the optimal solution.

**The type of the required elementary row transformations:** To make an entry in a column of $R_N$ nonnegative we require the following type of elementary row transformation: Suppose in r-th column the entry $b_{i_m}$ is positive and suppose the entry $b_{i_n}$ in that column is negative then we replace the row containing this negative entry, $b_{i_n}$, by the row made up of the addition of this same row containing this negative entry with the row obtained by multiplying the above mentioned row containing positive entry, $b_{i_m}$ , by the ratio [-( $b_{i_n}/b_{i_m}$ )]. This action changes this negative entry, $b_{i_n}$ , to zero.

One can now easily see that by proceeding on exactly similar lines one can develop the formulation as above also for minimization problems and can deal with them successfully.

## The Main Result:

We now state and prove the main result of this paper.

**Theorem 2.1:** The complexity of the above given algorithm (**Algorithm 2.1**) for solving a linear programming problem is of polynomial order and it is same as that of the complexity required for finding reduced-row-echelon-form.

**Proof:** The algorithm basically finds the reduced-row-echelon-form (rref) of the augmented matrix defining the problem, i.e., $R$ = rref (**[E, F]**). Now, whether the linear programming problem at hand is either inconsistent, or unbounded, or infeasible, or a valid one with a solution becomes clear from the nature of $R$, as mentioned in the algorithm.

When the linear programming problem at hand is either inconsistent, or unbounded, or infeasible nothing is to be done and the complexity in this case is clearly that of the complexity for the algorithm to find rref.

When the linear programming problem at hand is a valid one with a solution then to find that solution we further require to carry out certain elementary row transformations to make all the negative entries in the columns corresponding to non-basic variables nonnegative. Now, to make all the negative entries in a row belonging to columns for non-basic variables one requires one elementary row transformation. Since there are in all "m" rows in ***R*,** one will require to perform at most "m" elementary row transformations, so the overall complexity still remains dominated by the order of the complexity required for rref algorithm.

### 3. Examples:

We now proceed with some examples for maximization problems:

**Example 2.1:** Maximize: $x+y$

Subject to: $x+2y\le 4$

$$-x+y\le 1$$

$$4x+2y\le 12$$

$$x,y\ge 0$$

**Solution:** For this problem we have

```
R = [   1,    0,    0,    0,  1/2,  -d+6    ]
    [   0,    1,    0,    0, -1/2, -6+2*d  ]
    [   0,    0,    1,    0,  1/2, 10-3*d   ]
    [   0,    0,    0,    1,    1, 13-3*d   ]
```

So, clearly,

$$R_N = \begin{bmatrix} 1, & 0, & 0, & 0, & 1/2, & -d+6 \\ 0, & 0, & 1, & 0, & 1/2, & 10-3*d \\ 0, & 0, & 0, & 1, & 1, & 13-3*d \end{bmatrix}$$

For this example the column forming coefficients for non-basic variable $s_3$ contains nonegative numbers. So, we set $s_3 = 0$. Clearly, $d_{\min}$ = 3.3333 = Optimal value for the objective function. Using this value of optimum we have $x = 2.66, y = 0.66$, $s_1 = 0, s_2 = 3, s_3 = 0$.

**Example 2.2:** We first consider the duel of the example suggested by E. M. L. Beale [3], which brings into existence the problem of **cycling** for the simplex method, and provide a solution as per the above new method which offers it directly without any cycling phenomenon.

Maximize: $0.75\, x_1 - 20\, x_2 + 0.5\, x_3 - 6\, x_4$

Subject to: $0.25\, x_1 - 8\, x_2 - x_3 + 9\, x_4 \leq 0$

$0.5\, x_1 - 12\, x_2 - 0.5\, x_3 + 3\, x_4 \leq 0$

$x_3 \leq 0$

$x_1,\ x_2,\ x_3,\ x_4 \geq 0$

**Solution:** For this problem we have the following

$$R = \begin{bmatrix} 1, & 0, & 0, & 0, & -22/3, & 38/3, & 4/3, & (-14/3)d+4/3 \\ 0, & 1, & 0, & 0, & -7/24, & 11/24, & 1/24, & (-5/24)d+1/2 \\ 0, & 0, & 1, & 0, & 0, & 0, & 1, & 1 \\ 0, & 0, & 0, & 1, & 1/18, & 1/18, & 1/9, & 1/9-(1/18)d \end{bmatrix}$$

So, clearly,

$$R_N = \begin{bmatrix} 1, & 0, & 0, & 0, & -22/3, & 38/3, & 4/3, & (-14/3)d+4/3 \\ 0, & 1, & 0, & 0, & -7/24, & 11/24, & 1/24, & (-5/24)d+1/24 \\ 0, & 0, & 0, & 1, & 1/18, & 1/18, & 1/9, & 1/9-(1/18)d \end{bmatrix}$$

We perform following elementary row transformations on *R*: Let us denote the successive rows of *R* by R(1), R(2), R(3), R(4). We change

(i) R(1) → R(1) + 132*R(4)
(ii) R(2) → R(2) + (126/24)*R(4)

This leads to new transformed *R* as follows:

$$
R = \begin{array}{llllllll}
[1, & 0, & 0, & 132, & 0, & 20, & 1\ 6, & 16\text{-}12d] \\
[0, & 1, & 0, & 21/4, & 0, & 3/4, & 5/8, & 5/8\text{-}(1/2)d] \\
[0, & 0, & 1, & 0, & 0, & 0, & 1, & 1 \quad ] \\
[0, & 0, & 0, & 1, & 1/18, & 1/18, & 1/9, & 1/9-(1/18)d\ ]
\end{array}
$$

In the transformed *R* we have nonnegative columns for all non-basic variables, which are now those corresponding to $x_4$, $s_2, s_3$. So, by setting $x_4 = s_2 = s_3 = 0$ and setting all expressions of type $c_{i_r} d + e_{i_r} = 0$ in the last column we find $d_{\min} = \min\{ d_{i_r} \}$ = 1.25. Using this value in the last column of the newly obtained transformed *R* we have: $x_1$ = 1.0000, $x_2$ = 0, $x_3$ = 1, $x_4$ = 0, $s_1$ = 0.7500, $s_2$ = 0, $s_3$ = 0, and the maximum value of $d$ = 1.2500.

**Example 2.3:** We now consider an unbounded problem. The new method directly implies the unbounded nature of the problem through the positivity of the coefficients of $d$ in matrix $R$ for the problem.

Maximize: $-x + 3y$

Subject to: $-x - y \le -2$

$$x - 2y \le 0$$

$$-2x + y \le 1$$

$$x, y \ge 0$$

**Solution:** The following is the matrix $R$:

$$
R = \begin{array}{llllll}
[ & 1, & 0, & 0, & 0, & -3/5,\ (1/5)d\text{-}3/5\ ] \\
[ & 0, & 1, & 0, & 0, & -1/5,\ (2/5)d\text{-}1/5\ ] \\
[ & 0, & 0, & 1, & 0, & -4/5,\ (3/5)d\text{-}14/5\ ] \\
[ & 0, & 0, & 0, & 1, & 1/5,\ 1/5+(3/5)d\ ]
\end{array}
$$

Here, all the coefficients of $d$ are positive. So, by setting variable $s_3$ = 0 we can see that we can assign any arbitrarily large value to variable $d$ without violation of nonnegativity constraints for variables. Thus, the problem has an unbounded solution.

**Example 2.4:** We now consider a problem having an infeasible starting basis. We see that new algorithm has no difficulty to deal with it.

Maximize: $3x + 2y$

Subject to: $x + y \le 4$

$$2x + y \le 5$$

$$x - 4y \le -2$$

$$x, y \ge 0$$

**Solution:** The following is the matrix $R$:

$R$ = [ 1, 0, 0, 0, 1/7, (2/7)d-2/7 ]
[ 0, 1, 0, 0, -3/14, 3/7+(1/14)d]
[ 0, 0, 1, 0, 1/14, 27/7-(5/14)d]
[ 0, 0, 0, 1, -1/14, 36/7-(9/14)d]

We perform

(i) R(4) → R(4) + R(3)
(ii) R(2) → R(2) + 3*R(3)

This leads to getting new transformed $R$ as follows:

$R$ = [ 1, 0, 0, 0, 1/7, (2/7)d-2/7 ]
[ 0, 1, 3, 0, 0, 12- d ]
[ 0, 0, 1, 0, 1/14, 27/7-(5/14)d]
[ 0, 0, 1, 1, 0, (63/7) – d ]

The first two columns of $R$ correspond to basic variables $x, y$. Since the columns corresponding to non-basic variables $s_1, s_3$ contain nonnegative entries in $R_N$ , so we set these variables to zero. From last row we have $s_2 = 0$ and $d_{\min} =$ 9. Also from first and second rows, $x = 1$, $y = 3$

**Example 2.5:** We now consider an infeasible problem.

Maximize: $3x + 2y$

Subject to: $2x - y \leq -1$

$-x + 2y \leq 0$

$x, y \geq 0$

**Solution:** The following is the matrix $R$:

$R$ = [ 1, 0, 0, -1/4, (1/4)d ]
[ 0, 1, 0, 3/8, (1/8)d ]
[ 0, 0, 1, 7/8, (-3/8)d-1 ]

Here, the coefficient of $d$ is negative only in the last row and so

$R_N$ = [ 0, 0, 1, 7/8, (-3/8)d-1].

We perform following elementary row transformations on *R*: Let us denote the successive rows of *R* by R(1), R(2), R(3), R(4). We change

(i) R(1) → (2/7)*R(4) + R(1)

This leads to

$R$ = [ 1, 0, 2/7, 0, 1/7*d-2/7]
[ 0, 1, 0, 3/8, (1/8)d ]
[ 0, 0, 1, 7/8, (-3/8)d-1]

Setting $S_1, S_2$ equal to zero, we have for consistency of last row $d = -(8/3)$ and using this value for $d$ we have $y = -(1/3)$. Thus, this problem is infeasible.

**Remark 2.1:** Klee and Minty [4], have constructed an example of a set of linear programs with $n$ variables for which simplex method requires $2^n - 1$ iterations to reach an optimal solution. Theoretic work of Borgwardt [5] and Smale [6] indicates that fortunately the occurrence of problems belonging to the class of Klee and Minty, which don't share the average behavior, is so rare as to be negligible. We now proceed to show that there is no problem of efficiency for new algorithm in dealing with the problems belonging to this class.

**Example 2.6:** We now consider a problem for which the simplex iterations are **exponential** function of the size of the problem. A problem belonging to the class described by Klee and Minty [4], containing $n$ variables requires $2^n - 1$ simplex steps. We see that the new method doesn't require any special effort

Maximize: $100x_1 + 10x_2 + x_3$

Subject to: $x_1 \leq 1$

$$20x_1 + x_2 \leq 100$$

$$200x_1 + 20x_2 + x_3 \leq 10000$$

$$x_1,\ x_2,\ x_3 \geq 0$$

**Solution:** The following as the matrix $R$:

$$R = \begin{bmatrix} 1, & 0, & 0, & 0, & 1/10, & -1/100, & -90+(1/100)d \\ 0, & 1, & 0, & 0, & -1, & 1/5, & 1900-(1/5)d \\ 0, & 0, & 1, & 0, & 0, & -1, & 2d-10000 \\ 0, & 0, & 0, & 1, & -1/10, & 1/100, & 91-(1/100)d \end{bmatrix}$$

We perform following elementary row transformations on *R*: Let us denote the successive rows of *R* by R(1), R(2), R(3), R(4). We change

(i) R(2) → 10R(1) +R(2), and
(ii) R(4) → R(1) +R(4)
(iii) R(1) → 10*R(1)
(iv) R(1) → R(1) + R(2)
(v) R(3) → R(3) +10R(2)

This changes R to

$$R = \begin{bmatrix} 10, & 1, & 0, & 0, & 1, & 0, & 100 \\ 0, & 1, & 0, & 0, & 0, & 1/10, & 1000-(1/10)*d \\ 0, & 10, & 1, & 0, & 0, & 0, & d \\ 1, & 0, & 0, & 1, & 0, & 0, & 1 \end{bmatrix}$$

which produces (after setting non-basic variables to zero as usual) the solution:

$x_1 = 0$, $x_2 = 0$, $x_3 = 10000$, $s_1 = 1$, $s_2 = 100$, $s_3 = 0$, and the maximum value of the objective function, $d = 10000$.

**Example 2.7:** Maximize: $4x + 3y$

Subject to: $x + 3.5y \le 9$

$2x + y \le 8$

$x + y \le 6$

$x, y \ge 0$

For this problem we get following $R$

$R =$ [ 1, 0, 0, 0, -3, d-18 ]
[ 0, 1, 0, 0, 4, 24-d ]
[ 0, 0, 1, 0, -11, -57+5/2*d]
[ 0, 0, 0, 1, 2, 20-d ]

We perform following elementary row transformations:

R(1) → R(1) + (3/2)R(4)
R(3) → R(3) + (11/2)R(4)

this leads to:

$R =$ [ 1, 0, 0, 0, 0, 12 – 1/2d]
[ 0, 1, 0, 0, 4, 24-d ]
[ 0, 0, 1, 11/2, 0, 53-3d ]
[ 0, 0, 0, 1, 2, 20-d ]

which leads to maximal basic feasible solution: $d = 17.667$,

$x = 3.1665, y = 1.667, s_1 = 0, s_2 = 1.1665, s_3 = 0$

**Example 2.7:** Maximize: 3x1+2x2+3x3+4x4+x5

Subject to:

4x1+3x2-2x3+2x4-x5 ≤ 12
2x1+3x2+x3+3x4+x5 ≤ 15
3x1+2x2+x3+2x4+5x5 ≤ 20
2x1+4x2+x3+6x4+x5 ≤ 25
x3 ≤ 3

A =
[ 3, 2, 3, 4, 1, 0, 0, 0, 0, 0, d ]
[ 4, 3, -2, 2, -1, 1, 0, 0, 0, 0, 12]
[ 2, 3, 1, 3, 1, 0, 1, 0, 0, 0, 15]
[ 3, 2, 1, 2, 5, 0, 0, 1, 0, 0, 20]
[ 2, 4, 1, 6, 1, 0, 0, 0, 1, 0, 25]
[ 0, 0, 1, 0, 0, 0, 0, 0, 0, 1, 3 ]

R = rref(A)

```
R =
[1,    0,    0,    0,    0,    0,    10/27,    -2/27,    -13/27,    -34/27,  13/27*d-317/27]
[0,    1,    0,    0,    0,    0,     7/9,     -1/18,     -1/9,       5/9,    85/9-7/18*d  ]
[0,    0,    1,    0,    0,    0,      0,        0,         0,         1,          3        ]
[0,    0,    0,    1,    0,    0,   -16/27,    1/54,     10/27,     -5/27,   5/27+7/54*d   ]
[0,    0,    0,    0,    1,    0,    -8/27,    7/27,      5/27,     11/27,  178/27-5/27*d  ]
[ 0,   0,    0,    0,    0,    1,   -79/27,   37/54,     46/27,    166/27, 1157/27-65/54*d]
```

and its sub-matrix $R_N$ is

$R_N$ =

```
[0,    1,    0,    0,    0,    0,     7/9,     -1/18,     -1/9,       5/9,    85/9-7/18*d  ]
[0,    0,    0,    0,    1,    0,    -8/27,    7/27,      5/27,     11/27,  178/27-5/27*d  ]
[0,    0,    0,    0,    0,    1,   -79/27,   37/54,     46/27,    166/27, 1157/27-65/54*d]
```

With the appropriately chosen three elementary row transformations this $R_N$ becomes new

$R_N$ =

```
[0,  35/27,  0,    0,   7/9,   0,    7/9,    7/54,       0,       28/27, 469/27-35/54*d ]
[0,   8/21,  0,    0,    1,    0,     0,     5/21,      1/7,      13/21,   214/21-1/3*d ]
[0,  79/21,  0,    0,    0,    1,     0,    10/21,      9/7,     173/21, 1646/21-8/3*d  ]
```

After appropriately re-embedding the rows of this new $R_N$ in original place in R we get

```
R =
[1,    0,    0,    0,    0,    0,    10/27,    -2/27,    -13/27,    -34/27,  13/27*d-317/27]
[0,  35/27,  0,    0,   7/9,   0,     7/9,     7/54,        0,       28/27, 469/27-35/54*d ]
[0,    0,    1,    0,    0,    0,      0,       0,          0,         1,          3        ]
[0,    0,    0,    1,    0,    0,   -16/27,   1/54,      10/27,     -5/27,   5/27+7/54*d   ]
[0,   8/21,  0,    0,    1,    0,      0,     5/21,       1/7,      13/21,   214/21-1/3*d ]
[0,  79/21,  0,    0,    0,    1,      0,    10/21,       9/7,     173/21,  1646/21-8/3*d  ]
```

Further two more elementary row transformations leads us to

```
R =
[1,  79/21,  0,    0,    0,   1,  10/27,  76/189, 152/189, 1319/189, -59/27*d+12595/189]
[0,  35/27,  0,    0,   7/9,  0,   7/9,    7/54,     0,      28/27, 469/27-35/54*d       ]
[0,    0,    1,    0,    0,   0,    0,       0,      0,        1,          3              ]
[0,  35/27,  0,    1,   7/9,  0,   5/27,    4/27,  10/27,    23/27, 158/9-14/27*d        ]
[0,   8/21,  0,    0,    1,   0,    0,      5/21,   1/7,     13/21,   214/21-1/3*d       ]
[0,  79/21,  0,    0,    0,   1,    0,     10/21,   9/7,    173/21,   1646/21-8/3*d      ]
```

Using this R and setting all non-basic variables to zero we can easily find the maximal value of the objective function, namely, $d = 26.8$ and by further substitution of this value of objective function in we can easily find the optimal solution.

4. **A New Algorithm for Nonlinear Programming:** We now proceed show that we can deal with nonlinear constrained optimization problems using the same above given technique used to deal with linear programming problems. The algorithms developed by Bruno Buchberger which transformed the abstract notion of Grobner basis into a fundamental tool in computational algebra will be utilized. The technique of Grobner bases is essentially a version of reduced row echelon form (used above to handle the linear programs made up of linear polynomials) for higher degree polynomials [7]. A typical nonlinear program can be stated as follows:

   Maximize/Minimize: $f(x)$

   Subject to: $h_j(x) = 0, j = 1,2,\cdots,m$

   $g_j(x) \geq 0, j = m+1, m+2,\cdots, p$

   $x_k \geq 0, k = 1,2,\cdots,n$

Given a nonlinear optimization problem we first construct the following nonlinear system of equations:

$$f(x) - d = 0 \tag{3.1}$$

$$h_j(x) = 0, j = 1,2,\cdots,m \tag{3.2}$$

$$g_j(x) + s_j = 0, j = m+1, m+2,\cdots, p \tag{3.3}$$

where $d$ is the unknown parameter whose optimal value is to be determined subject to nonnegativity conditions on problem variables and slack variables. For this to achieve we first transform the system of equations into an equivalent system of equations bearing the same solution set such that the system is easier to solve. We have seen so far that the effective way to deal with linear programs is to obtain the reduced row echelon form for the combined system of equations incorporating objective equation and constraint equations. We will see that for the nonlinear case the effective way to deal with is to obtain the equivalent of reduced row echelon form, namely, the Grobner basis representation for this system of equations (3.1)-(3.3). We then set up the equations obtained by equating the partial derivatives of $d$ with respect to problem variables $x_i$ and slack variables $s_i$ to zero and utilize the standard theory and methods used in calculus. We demonstrate the essence of this method by solving certain examples. where $d$ is the unknown parameter whose optimal value is to be determined subject to nonnegativity conditions on problem variables and slack variables. For this to achieve we first transform the system of equations into an equivalent system of equations bearing the same solution set such that the system is easier to solve. We have seen so far that the effective way to deal with linear programs is to obtain the reduced row echelon form for the combined system of equations incorporating objective equation and constraint equations. We will see that for the nonlinear case the effective way to deal with is to obtain the equivalent of reduced row echelon form for the set of polynomials, namely, the Grobner basis representation for this system of equations (3.1)-(3.3). We then set up the equations obtained by equating the partial derivatives of $d$ with respect to problem variables $x_i$ and slack variables $s_i$ to zero and utilize the standard theory and methods used in calculus. We demonstrate the essence of this method by solving an example:
These examples are taken from [8], [9]. These examples sufficiently illustrate the power of this new method of using powerful technique of Grobner basis to successfully and efficiently deal with nonlinear programming problems.

**Example 3.1:** Maximize: $-x_1^2+4x_1+2x_2$

Subject to: $x_1+x_2\le 4$

$2x_1+x_2\le 5$

$-x_1+4x_2\ge 2$

**Solution:** We build the following system of equations:

$$-x_1^2+4x_1+2x_2-d=0$$
$$x_1+x_2+s_1-4=0$$
$$2x_1+x_2+s_2-5=0$$
$$-x_1+4x_2-s_3-2=0$$

such that: $x_1,x_2,s_1,s_2,s_3\ge 0$

We now transform the nonlinear/linear polynomials on the left hand side of the above equations by obtaining Grobner basis for them as follows:

$$486-81d-18s_2-16s_2^2+36s_3-8s_2s_3-s_3^2=0 \quad (3.1.1)$$
$$9-9s_1+5s_2-s_3=0 \quad (3.1.2)$$
$$-9+s_2-2s_3+9x_2=0 \quad (3.1.3)$$
$$-18+4s_2+s_3+9x_1=0 \quad (3.1.4)$$

In order to maximize the objective function represented in terms of parameter d the partial derivatives of d with respect to problem variables must vanish, i.e. they must be equal to zero. therefore,

By setting $\frac{\partial d}{\partial s_2}=0$ and $\frac{\partial d}{\partial s_3}=0$ we get equations:

$$32s_2+8s_3=-18$$
$$8s_2+2s_3=36$$

a rank deficient system. Note that for maximization of $d$ if we set $\frac{\partial d}{\partial s_2}=0$ we get the value of $s_2$ that maximizes $d$, namely, $s_2=-(18/32)-(8s_3/32)$, a negative value for any nonnegative value of $s_3$. So, we set $s_2=0$. Similarly, for maximization of $d$ if we set $\frac{\partial d}{\partial s_3} = 0$ we get the value of $s_3$ that maximizes $d$, namely, $s_3 = 18 - 4s_2 (= 18)$, setting $s_2 = 0$. But, by setting $s_2 = 0$ in the second equation above the largest possible value for $s_3$ that one can have (is obtained by setting $s_1 = 0$and it) is **9**, when $s_2 = 0$. Thus, setting $s_2 = 0, s_3 = 9$ in the first equation we get $d$ = **9**. From third and fourth equation we get $x_2 = 3, x_1 = 1$.

**Example 3.2:** Maximize: $-8x_1^2 - 16x_2^2 + 24x_1 + 56x_2$

Subject to: $x_1 + x_2 \leq 4$

$$2x_1 + x_2 \leq 5$$

$$-x_1 + 4x_2 \geq 2$$

$$x_1, x_2 \geq 0$$

**Solution:** We build the following system of equations:

$$-8x_1^2 - 16x_2^2 + 24x_1 + 56x_2 - d = 0$$

$$x_1 + x_2 + s_1 - 4 = 0$$

$$2x_1 + x_2 + s_2 - 5 = 0$$

$$-x_1 + 4x_2 - s_3 - 2 = 0$$

We now transform the nonlinear/linear polynomials on the left hand side of the above equations by obtaining Grobner basis for them as follows:

$$504 - 9d + 8s_2 - 16s_2^2 + 56s_3 - 8s_3^2 = 0 \qquad (3.2.1)$$

$$9 - 9s_1 + 5s_2 - s_3 = 0 \qquad (3.2.2)$$

$$-9 + s_2 - 2s_3 + 9x_2 = 0 \qquad (3.2.3)$$

$$-18 + 4s_2 + s_3 + 9x_1 = 0 \qquad (3.2.4)$$

from first equation (3.2.1), in order to maximize $d$, we determine the values of $s_2, s_3$ as follows:

If we set $\frac{\partial d}{\partial s_2} = 0$ we get the value of $s_2$ that maximizes $d$, namely, $s_2 = \frac{1}{4}$.

Similarly, if we set $\frac{\partial d}{\partial s_3} = 0$ we get the value of $s_3$ that maximizes $d$, namely, $s_3 = \frac{7}{2}$. Putting these values of $s_2, s_3$ in the first and second equation we get respectively the maximum value of $d = 67$ and the value of $s_1 = \frac{3}{4}$. Using further these values in the third and fourth equation we get $x_1 = 1.5, x_2 = 1.75$.

**Example 3.3:** Minimize: $(x_1 - 3)^2 + (x_2 - 4)^2$

Subject to: $2x_1 + x_2 = 3$

**Solution:** We form the objective equation and constraint equations as is done in the above examples and then find the Grobner basis which yields:

$$5x_1^2 - d - 2x_1 + 10 = 0$$
$$2x_1 + x_2 - 3 = 0$$

Setting $\frac{\partial d}{\partial x_1} = 0$ we get the value of $x_1$ that minimizes $d$, namely, $x_1 = 0.2$.

This yields $d = 9.8$ and $x_2 = 2.6$

**Example 3.4:** Minimize: $x_1^2 - x_2$

Subject to: $x_1 + x_2 = 6$

$x_1 \geq 1$

$x_1^2 + x_2^2 \leq 26$

**Solution:** We form the objective equation and constraint equations as is done in the above examples and then find the Grobner basis which yields:

$$2d - 14s_1 - s_2 + 8 = 0$$
$$x_2 + s_1 - 5 = 0$$
$$x_1 - s_1 - 1 = 0$$
$$-s_2 - 8s_1 + 2s_1^2 = 0$$

For minimizing $d$ we should set the values of $s_1, s_2$ equal to zero (as they have signs opposite to $d$) which yields $d = -4$. From other equations we get $x_1 = 1, x_2 = 5$.

**Example 3.5:** Minimize: $-6x_1 + 2x_1^2 - 2x_1x_2 + 2x_2^2$

Subject to: $x_1 + x_2 \leq 2$

**Solution:** We form the objective equation and constraint equations as is done in the above examples and then find the Grobner basis which yields:

$$6x_2^2 - d + 6x_2s_1 - 6x_2 + 2s_1^2 - 2s_1 - 4 = 0 \quad (3.5.1)$$
$$x_1 + x_2 + s_1 - 2 = 0 \quad (3.5.2)$$

Setting $\frac{\partial d}{\partial x_2} = 0$ and $\frac{\partial d}{\partial s_1} = 0$ we get equations:

$$12x_2 + 6s_1 = 6$$
$$6x_2 + 4s_1 = 2$$

There solution gives $s_1 = -1$, which is forbidden so we first set $s_1 = 0$ in the initial equations (3.5.1) and (3.5.2) and again set $\frac{\partial d}{\partial x_2} = 0$ which yields the value of $x_2$ that minimizes $d$, namely, $x_2 = \frac{1}{2}$. This in turn produce $x_1 = \frac{3}{2}$ and $d = -\frac{11}{2}$.

**5. Conclusion:** Condensing of the linear form (to be optimized) into a new parameter and developing the appropriate equations containing it is a useful idea. This idea is useful not only for linear programs but also for nonlinear programming problems and provides new effective ways to deal with these problems.

## Acknowledgements

The author is thankful to Dr. M. R. Modak for useful discussion.